\documentclass[11pt]{article}
\usepackage{amsmath,amsthm,amsfonts,amssymb,mathrsfs,epsfig,bm}
\usepackage[usenames]{color}

\parskip 	\smallskipamount

\newcommand{\keywords}[1]{ \noindent {\footnotesize
             {\small \em Keywords and phrases.} {\sc #1} } }
\newcommand{\ams}[2]{  \noindent {\footnotesize
             {\small \em AMS {\rm 2000} subject classifications.
             {\rm Primary {\sc #1}; secondary {\sc #2}} } } }

\newtheorem{theorem}{Theorem}
\newtheorem{lemma}{Lemma}
\newtheorem{corollary}{Corollary}

\newtheorem{definition}{Definition}
\theoremstyle{remark}

\def\N{\mathbb{N}}
\def\Z{\mathbb{Z}}

\def\R{\mathbb{R}}

\def\AA{\mathscr{A}}

\def\EE{\mathscr{E}}

\def\NN{\mathcal{N}}
\def\PP{\mathscr{P}}
\renewcommand{\phi}{\varphi}
\renewcommand{\epsilon}{\varepsilon}

\newcommand{\1}{{\text{\Large $\mathfrak 1$}}}

\renewcommand{\emptyset}{\varnothing}
\newcommand{\I}{_{\text{\sc i}}}
\newcommand{\II}{_{\text{\sc ii}}}
\newcommand{\III}{_{\text{\sc iii}}}
\newcommand{\IV}{_{\text{\sc iv}}}

\def\nn{\text{\tiny $N$}}
\newcommand{\bt}[6]{{\tt BITTORRENT[$#1,#2,#3,#4,#5,#6$]}}
\newcommand{\ttt}{\widetilde}

\begin{document}

\title{\bf A stochastic epidemiological model and a deterministic
limit for BitTorrent-like peer-to-peer file-sharing networks}
\author{
{\sc George Kesidis}
\thanks{CSE and EE Depts, Pennsylvania State University,
University Park, PA 16802, USA;~~ {\tt gik2@psu.edu}
}\\
Pennsylvania State University
\and
{\sc Takis Konstantopoulos}
\thanks{School of Mathematical \& Computer Sciences and the Maxwell
Institute for Mathematical Sciences,
Heriot-Watt University, Edinburgh EH14 4AS, UK;~~ {\tt takis@ma.hw.ac.uk}
}\\
Heriot-Watt University
\and
{\sc Perla Sousi}
\thanks{Statistical Laboratory, Centre for Mathematical Sciences,
University of Cambridge,
Wilberforce Road,
Cambridge CB3 0WB, UK;~~
{\tt P.Sousi@statslab.cam.ac.uk}
}\\
University of Cambridge
}
\date{}
\maketitle

\begin{abstract}
In this paper, we propose a stochastic model for a
file-sharing peer-to-peer network
which resembles the popular BitTorrent system:
large files are split into chunks and a peer can download
or swap from another peer only one chunk at a time. 
We prove that the fluid limits of a scaled Markov model of this system
are of the coagulation form, special cases of which are 
well-known epidemiological (SIR) models.
In addition, Lyapunov stability and settling-time results are
explored. We derive conditions under which the
BitTorrent incentives under consideration result in
shorter mean file-acquisition times for peers compared
to client-server (single chunk) systems.
Finally, a diffusion approximation is given and some open
questions are discussed.

\vspace*{2mm}
\keywords{fluid limit, peer-to-peer networks, diffusion approximation, stability, epidemiological model}

\vspace*{2mm}
\ams{60F17,90B15}{60J75}

\end{abstract}

\section{Introduction}
\renewcommand{\thefootnote}{\fnsymbol{footnote}}

Peer-to-peer (p2p) activity continues to represent a very significant fraction
of  overall Internet traffic, 44\% by one recent account 
\cite{DCInfo063008}.
BitTorrent
\cite{BitTorrent,Cohen:2003,Ge:2004,Yang:2004,Qiu:2004,bittor-imc05,Turner:2005}
is a widely deployed p2p
file-sharing network which has recently played a significant
role in the network neutrality debate.
Under BitTorrent, peers join ``swarms" (or ``torrents")
where each swarm corresponds
to a specific data object (file). 
The process of finding the peers 
in a given swarm to connect to 
is typically facilitated through a centralised ``tracker".
Recently, a trackerless BitTorrent client
has been introduced that uses distributed hashing for query
resolution \cite{maymounkov02kademlia}.  

For file sharing, a peer is typically uploads
upload pieces (``chunks") of the file to other peers in the
swarm while downloading his/her missing chunks from them. 
This chunk swapping constitutes a transaction-by-transaction
incentive for peers to cooperate (i.e., trading rather
than simply download) to disseminate data objects. 
Large files may be segmented 
into several hundred chunks, all of which the  peers of the
corresponding warm must collect,
and in the process disseminate their own chunks
before they can reconstitute the desired file and possibly leave
the file's swarm.

In addition to the framework
in which data objects are segmented into chunks to promote 
cooperation through swapping,  there is a system whereby the rate
at which chunks are
uploaded is assessed for any given transaction, and peers that
allocate inadequate bandwidth for uploading may be ``choked"
\cite{LZhang07,ITA06}. Choking may also be applied to peers 
who, by employing multiple identities (sybils), 
abuse BitTorrent's  system of allowing newly arrived peers to
a swarm  to just download a few chunks (as they clearly cannot
trade what they simply do not as yet possess).
BitTorrent can also rehabilitate peers by (optimistically) unchoking them.
In the following, we do not directly consider upload bandwidth and
related choking issues.

In this paper, we motivate a deterministic epidemiological model
of file dissemination for peer-to-peer file-sharing networks that
employ BitTorrent-like incentives, a generalisation of that
given in \cite{Kesidis:2006b}\footnote{And this paper is a significant
extension of \cite{ICASSP07}.}.
Our model is different from those explored in 
\cite{Mass-Voj-05,Yang:2004,Qiu:2004} for BitTorrent,
and we compute different quantities of interest. Our epidemiological 
framework, similar to that we used for the spread of multi-stage worms
\cite{Kesidis08},  could also be adapted
for network coding systems.
In \cite{bittor-imc05}, the authors propose a ``fluid" model 
of a single torrent/swarm (as we do in the following) and fit it
to (transient) data drawn from aggregate swarms. 
The connection to branching
process models \cite{Yang:2004,Ge:2004} is simply
that ours only tracks the number of active peers who possess or
demand the file under consideration, i.e., a single swarm. 
Though our model is
significantly simpler than that of prior work, 
it is derived directly from an intuitive transaction-by-transaction
Markov process modelling file-dissemination of the p2p network and its
numerical solutions clearly demonstrate the effectiveness of the
aforementioned incentives. A basic assumption in the following is
that peers do not distribute bogus files (or file chunks)
\cite{Walsh:2005:2}.

This paper is organised as follows. 
The Markovian model is developed in detail in Section \ref{stoch-sec}.
A proof of its fluid limit is described in 
Section \ref{fluid-limit-sec}  (and in the Appendix)
including Lyapunov stability and settling time results. 
The behaviour of the limiting ODE is studied in Section 
\ref{ODEbeh}  for specific examples.
In Section \ref{incentives-sec}, we derive conditions under
which the BitTorrent-like system
has improved performance (smaller mean time to completely
acquire the file by a peer) compared to a system of
pure client-server (no chunk-swapping) interactions.
A diffusion approximation is given in \ref{diff-sec}.
The paper concludes with a discussion of open problems in 
Section \ref{concl-sec}.

\section{The stochastic model}\label{stoch-sec}
We fix a set $F$ (a file) which is partitioned into $n$ 
(on the order of hundreds) pieces called chunks.
Consider a large networked ``swarm" of $N$ nodes called peers. 
Each peer possesses a certain (possibly empty) subset $A$ of $F$.
As time goes by, this peer interacts with other peers,
the goal being to enlarge his set $A$ until, eventually,
the peer manages to collect all $n$ chunks of $F$.
The interaction between peers can either be a download
or a swap; in both cases, chunks are being copied from
peer to peer and are assumed never lost.
A peer will stay in the network as long as he does not possess all chunks.
After collecting everything, sooner or later a peer departs
or switches off.
By splitting the desired file into many chunks we give incentives
to the peers to remain active in the
swarm for long time during which other peers will
take advantage of their possessions. 

\subsection{Possible interactions}
We here describe how two peers, labelled $A, B$, interact.
The following types of interactions are possible:
\begin{enumerate}
\item
{\bf Download:} Peer $A$ downloads a chunk $i$ from $B$. This is possible
only if $A$ is a strict subset of $B$. 
If $i \in B$ then, after the downloading
$A$ becomes $A' = A \cup \{i\}$ and but $B$ remains $B$ because it
since it gains nothing from $A$.
Denote this interaction as:
\[ 
\boxed{
(A \leftarrow B) \rightsquigarrow (A', B)
}
\]
The symbol on the left is supposed to show the type
of interaction and the labels before it,
while the symbol on the right shows the labels after the interaction.
\item
{\bf Swap:} Peer $A$ swaps with peer $B$. In other words, $A$ gets a chunk
$j$ from $B$ and $B$ gets a chunk $i$ from $A$. 
It is required that $j$ is not an element of $A$ and $i$ not 
an element of $B$.
We denote this interaction by
\[
\boxed{
(A \leftrightarrows B) \rightsquigarrow (A', B')
}
\]
where $A' = A \cup \{j\}$, $B' = B \cup \{i\}$.
We thus need $A \setminus B \neq \varnothing$ and $B \setminus A \neq
\varnothing$.
\end{enumerate}

\subsection{Notation}
The set of all combinations of $n$ chunks, which partition
$F$, is denoted by $\PP(F)$, where $|\PP(F)|=2^n$ and the empty set is
included.
We write $A \subset B$ (respectively, $A \subsetneq B$) when $A$
is a subset (respectively, strict subset) of $B$.
We (unconventionally) write
\begin{equation*}
\text{ $A \sqsubset A'$ when $A \subset A'$ and $|A'-A|=1$. }
\end{equation*}
If $A \cap B = \varnothing$, we use $A + B$ instead of $A \cup B$;
if $B=\{b\}$ is a singleton, we often write $A+b$ instead of $A + \{b\}$.
If $A \subset B$ we use $B-A$ instead of $B \setminus A$.
We say that \begin{equation*}
\text{ $A$ relates to $B$ (and write $A \sim B$)
 when $A \subset B$ or $B \subset A$};
\end{equation*}
if this is not the case, we write $A \not \sim B$.
Note that $A \not \sim B$ if and only if two peers labelled $A$, $B$ 
can swap chunks.
The space 
of functions (vectors)
from $\PP(F)$ into $\Z_+$ is denoted by $\Z_+^{\PP(F)}$.
The stochastic model will take values in this space.
The deterministic model will evolve in $\R_+^{\PP(F)}$.
We let $e_A \in \Z_+^{\PP(F)}$ be the vector with
coordinates \[
e_A^B := \1(A=B), \quad B \in \PP(F).
\]
For $x \in \Z_+^{\PP(F)}$ or $\R_+^{\PP(F)}$
we let $|x| := \sum_{A \in \PP(F)} |x^A|$.
If $\AA \subset \PP(F)$ then the $\AA$-face $\R_+^\AA$ of $\R_+^{\PP(F)}$ is
defined by $\R_+^\AA:=\{x \in \R_+^{\PP(F)}:~ 
\sum_{A \in \AA} x^A=0,~ \prod_{B \not \in \AA} x^B >0\}$.

\subsection{Defining the rates of individual interactions}
We follow the logic of stochastic modelling of chemical reactions
or epidemics and assume that the chance of a particular interaction
occurring in a short interval of time is proportional to the number
of ways of selecting the peers needed for this interaction
\cite{Kurtz81}.
Accordingly, the interaction rates 
{\em must} be given by the formulae described below.

Consider first finding the rate of
a download $A \leftarrow B$, where $A \subsetneq B$,
when the state of the system is $x
 \in \Z_+^{\PP(F)}$.
There are $x^A$ peers labelled $A$ and $x^B$ labelled $B$. 
We can choose them in $x^A x^B$ ways.
Thus the rate of a download $A \leftarrow B$
that results into $A$ getting {\em some} chunk from $B$ should be
proportional to $x^A x^B$. 
However, we are interested in the rate
of the {\em specific} 
interaction $(A \leftarrow B) \rightsquigarrow (A', B)$,
that turns $A$ into a specific set $A'$ differing from $A$
by one single chunk $(A \sqsubset A')$; 
there are $|B -A|$ chunks that $A$ can download from $B$;
the chance that picking one of them is $1/|B-A|$. 
Thus we have:
\[
(DR) \quad
\left\{
\text{
\begin{minipage}{12cm}
\em
the rate of the download $(A \leftarrow B) \rightsquigarrow (A', B)$
equals $\displaystyle \beta x^A \frac{x^B}{|B-A|}$,
\\
as long as 
$A \sqsubset A' \subset B$, 
\end{minipage}
}
\right.
\]
where $\beta > 0$.

Consider next a swap $A \leftrightarrows B$ 
and assume the state is $x$. 
Picking two peers labelled $A$ and $B$ (provided that $A \not \sim B$)
from the population is done in $x^A x^B$ ways. 
Thus the rate of a swap $A \leftrightarrows B$
is proportional to $x^A x^B$.
So if we {\em fix} two chunks
$i \in A\setminus B, j \in B\setminus A$ 
and specify that $A'=A+j, B'=B+i$, then 
the chance of picking $i$ from $A\setminus B$ and $j$ from $B\setminus A$
is $1/|A\setminus B| |B\setminus A|$.
Thus,
\[
(SR) \quad \left\{
\text{
\begin{minipage}{12cm}
\em
the rate of the swap $(A \leftrightarrows B) \rightsquigarrow (A', B')$
equals $\displaystyle \gamma 
\frac{x^A x^B}{|A \setminus B| |B \setminus A|}$,
\\
a long as 
$A \sqsubset A',
\quad
B \sqsubset B',
\quad
A'-A \subset B,
\quad
B'-B \subset A$,
\end{minipage}
}
\right.
\]
where $\gamma > 0$.

%

\subsection{Deriving the Markov chain rates}
Having defined the rates of each individual interaction 
we can easily define rates
$q(x,y)$
of a Markov chain in continuous time and state space $\Z_+^{\PP(F)}$
as follows. 

Define functions $\lambda_{A,A'}, \mu_{A,B} :
\R^{\PP(F)} \to \R$ by:
\begin{subequations}		\label{lllmmmm}
\begin{align}			
&\lambda_{A,A'}(x) 
:=
\left[\beta x^A \sum_{C : C \supset A'}
\frac{x^C}{|C-A|} \right]
~\1(A \sqsubset A') 			\label{lll}
\\
&\mu_{A,B}(x)		
:= \gamma 
\frac{x^A x^B}{|A \setminus B| |B \setminus A|}
~\1(A \not \sim B). 			\label{mmm}
\end{align}
\end{subequations}
Consider also constants $\delta \ge 0$ and $\alpha^A \ge 0$ for $A \in \PP(F)$,
i.e., $\alpha \in\R_+^{\PP(F)}$.
The transition rates of the closed conservative Markov chain are given by:
\begin{equation}			\label{rates}
q(x,y) :=
\begin{cases}
\lambda_{A,A'}(x), 
& \text{ if } y = x-e_A+e_{A'}
\\
\mu_{A,B}(x), 
& \text{ if } 
\begin{cases}
y = x-e_A-e_B+e_{A'}+e_{B'}
\\
A \sqsubset A', B \sqsubset B', A'-A \subset B, B'-B \subset A,
\end{cases}
\\
\alpha^A
& \text{ if } y = x+e_A
\\
\delta x^F
& \text{ if } y = x-e_F
\\
0, & \text{ for any other value of } y \not = x,
\end{cases}
\end{equation}
where $x$ ranges in $\Z_+^{\PP(F)}$.

A little justification of the first two cases is needed:
that $q(x,x-e_A-e_B+e_{A'}+e_{B'}) = \mu_{A,B}(x)$
is straightforward. It corresponds to a swap,
which is only possible when 
$A \sqsubset A', B \sqsubset B', A'-A \subset B, B'-B \subset A$. 
The swap rate was defined by (SR).
To see that $q(x,x-e_A+e_{A'})=\lambda_{A,A'}(x)$ we observe
that a peer labelled $A$ can change its label to $A' \sqsupset A$ 
by downloading a chunk from some
set $C$ that contains $A'$, so we sum the rates (DR) over all
these possible individual interactions to obtain the first line
in \eqref{rates}.
We can think of having Poisson process of arrivals of new peers 
at rate $|\alpha|$, and that each arriving peer is labelled $A$
with probability $\alpha^A/|\alpha|$.
Peers can depart, by definition, only when they are labelled $F$
and it takes an exponentially distributed amount of time (with mean
$1/\delta$) for a departure to occur. Thus, $q(x, x-e_F)=\delta x^F$.
We shall let $\mathsf Q$ denote the generator of the chain, i.e.\
$\mathsf Q f(x) = \sum_y (f(y)-f(x)) q(x,y)$, when $f$ is an appropriate functional
of the state space.

\begin{definition}
[\bt{x_0}{n} {\alpha} {\beta}{\gamma} {\delta}]
Given $x_0 \in \Z_+^{\PP(F)}$ (initial configuration),
$n=|F| \in \N$ (number of chunks),
$\alpha\in \R_+^{\PP(F)}$ (arrival rates), $\beta > 0$
(download rate), $\gamma \ge0$ (swap rate), $\delta \ge 0$ (departure rate)
we let \bt{x_0}{n} {\alpha} {\beta}{\gamma} {\delta}
be a Markov chain $(X_t, t \ge 0)$ with 
transition rates \eqref{rates} and $X_0=x_0$.
We say that the chain (network) is \underline{open} 
if $\alpha^A > 0$ for at least one $A$
and $\delta >0$;
it is \underline{closed} if $\alpha^A=0$ for all $A$;
it is \underline{conservative} 
if it is closed and $\delta=0$;
it is \underline{dissipative} if it is closed and $\delta>0$.
\end{definition}

In a conservative network, we have $q(x,y)=0$ if $|y| \neq |x|$ and
so $|X_t| = |X_0|$ for all $t \ge 0$.
Here, the actual state space is the simplex
\[
\{x \in \Z_+^{\PP(F)}: |x|=N\},
\]
where $N = |X_0|$.
It is easy to see that 
the state $e_F$ is reachable from any other state, but
all rates out of $e_F$ are zero. Hence a conservative network has
$e_F$ as a single absorbing state.

In a dissipative network, we have
$|X_t| \le |X_0|$ for all $t \ge 0$.
Here the state space is 
\[
\{x\in \Z_+^{\PP(F)}: |x|\le N\},
\]
where $N = |X_0|$.
It can be seen that a dissipative network has many absorbing points.

In an open network, there are no absorbing points. On the other
hand, one may wonder if certain components can escape to infinity.
This is not the case:
\begin{lemma}
If $\alpha^F > 0$ then
the open {\tt BITTORRENT[$x,n,\beta,\gamma,\alpha,\delta$]}
is positive recurrent Markov chain.
\end{lemma}
\proof (sketch)
If $\alpha^F > 0$, $\delta > 0$ the Markov chain is irreducible.
The remainder of the proof is based on a the construction of
a simple Lyapunov function:
\[
V(x) := |x|,
\]
for which it can be shown that there is a bounded set of states $K$ 
such that 
\[
\sup_{x \not \in K} (\mathsf Q V)(x) < 0.
\]
Perhaps the easiest way to
see this is by appealing to the stability of the corresponding ODE limit;
see Theorem \ref{ODEapprox} below and \cite{FK}.
\qed

\subsection{Example: $n=1$}
Let us take the special case where the file consists of a single chunk
$(n=1)$. The state here is $x=(x^\varnothing, x^1:=x^F)$. 
The rates are:
\begin{align}
&q\big((x^\varnothing, x^1), (x^\varnothing+1, x^1)\big)
= \alpha^\varnothing \nonumber \\
&q\big((x^\varnothing, x^1), (x^\varnothing, x^1+1)\big)
= \alpha^1 \nonumber \\
&q\big((x^\varnothing, x^1), (x^\varnothing-1, x^1+1)\big)=\beta x^\varnothing x^1  \nonumber \\
&q\big((x^\varnothing, x^1), (x^\varnothing, x^1-1)\big)=\delta x^1.
\label{kmc}
\end{align}
If $\alpha^\varnothing=\alpha^1=0$,
this is the stochastic version of the classical (closed) Kermack-McKendrick
(or susceptible-infective-removed (SIR)) 
model for a simple epidemic process \cite{Daley-Gani}.
Its absorbing points are states of the form $(x^\varnothing, 0)$.
In epidemiological terminology, $x^1$ is the number of infected
individuals, whereas $x^\varnothing$ is the number of susceptible
ones. Contrary to the epidemiological interpretation, infection
{\em is} desirable, for infection is tantamount to downloading
the file. 

\section{Macroscopic description: fluid limit}\label{fluid-limit-sec}
Analysing the Markov chain in its original form is complicated. 
We thus resort to a first-order approximation by an ordinary differential
equation (ODE). 

Let $v(x)$ be the vector field on $\R_+^{\PP(F)}$ with
components $v^A(x)$ defined by
\begin{multline} 
\label{vcomp}
\lefteqn{
v^A(x) = 
\alpha^A
-x^A \big(\beta \phi_d^A(x) + \gamma \phi^A_s(x) \big)
}
\\
+ \beta
\sum_{B: A \subset B}
\frac{\psi^A_d(x) x^B }{1+|B\setminus A|}
+
\gamma
\sum_{B: A \not \subset B}
\frac{\psi^{A,B}_s(x)  x^B }{1+|B\setminus A|}
-\delta x^F \1(A=F),
\end{multline}
where
\begin{align}
& \phi_d^A(x) := \sum_{B \supset A} x^B, \quad 
\phi_s^A(x) := \sum_{B \not \sim A} x^B 
\nonumber   \\
& \psi_d^A(x) := \sum_{a \in A} x^{A-a}, \quad 
\psi_s^{A,B}(x) := \sum_{a \in A \cap B} x^{A-a} 
\label{phipsi}
\end{align}
Consider the differential equation
\begin{equation}
\dot x = v(x) \text{ with initial condition $x_0$.} \label{ODE0}
\end{equation}
Consider the sequence of stochastic models \bt{X_{\nn,0}}{n}
{N \alpha}{\frac{\beta}{N}}{\frac{\gamma}{N}}{\delta} for
$N \in \N$ and let $X_{\nn,t}$ be the corresponding jump Markov chain.

\begin{theorem}\label{ODEapprox}
There is a has a unique smooth
(analytic) solution to (\ref{ODE0}), denoted by $x_t$ for $t \ge 0$.
Also, if there is an $x_0 \in \R_+^{\PP(F)}$ such that
$X_{\nn,0}/N \to x_0$ 
as $N \to \infty$,
then for any $T, \epsilon > 0$,
\[
\lim_{N \to \infty} 
P\big(\sup_{0 \le t \le T}|N^{-1} X_{\nn,t} -x_t| > \epsilon\big)=0.
\]
\end{theorem}

\proof
Let $\NN$ be the set of vectors $-e_{F}$, $e_A$, 
$-e_A+e_{A'}$, $-e_A-e_B+e_{A'}+e_{B'}$,
where $A, B\in\PP(F)$ and
$A \sqsubset A'$, $B \sqsubset B'$.
{From} \eqref{rates}, we have that $q(x,y)=0$ if $y-x \not \in \NN$.
Introduce, for each $\zeta \in \NN$, a unit rate Poisson process
$\Phi_\zeta$ on the real line, and assume that these Poisson processes
are independent.
Consider the Markov chain $(X_t)$ for the 
\bt{X_0}{n}{\alpha}{\beta}{\gamma}{\delta}.
Its rates are of the form
\begin{equation}
\label{Qdef}
q(x,x+\zeta) = Q_\zeta(x), \quad \zeta \in \NN,
\end{equation}
where $Q_\zeta(x)$ is a polynomial in $2^n$ variables of degree $2$, and
which can be read directly from \eqref{rates}; its coefficients
depend on the parameters $\alpha$, $\beta$, $\gamma$, $\delta$.
We can represent \cite{Kurtz81,Kurtz86} $(X_t)$ as:
\[
X_t=X_0 + \sum_{\zeta \in \NN} \zeta
\Phi_\zeta\big(\int_0^t Q_\zeta(X_s) ds\big).
\]
Consider now the Markov chain $\frac{1}{N} X_{\nn,t}$ corresponding to
to \bt{X_{\nn,0}}{n}{(N \alpha^{A})}{\beta/N}
{\gamma/N}{\delta}.
The transition rates for $\frac{1}{N} X_{\nn,t}$ are
\[
q(x/N, (x+\zeta)/N) = N Q_\zeta(x/N), \quad
x \in \Z_+^{\PP(F)}, \quad \zeta \in \NN,
\]
and $0$, otherwise.
Here, $Q_\zeta(x)$ is the polynomial defined through 
\eqref{Qdef} and \eqref{rates} and we now assume that its variables
range over the reals.
Therefore, $\frac{1}{N} X_{\nn,t}$ can be represented as
\[
\frac{1}{N}X_{\nn,t}=\frac{1}{N}X_{\nn,0} + \sum_{\zeta \in \NN} \zeta
\frac{1}{N} 
\Phi_\zeta\bigg(N \int_0^t Q_\zeta(\frac{1}{N}X_{\nn,s}) ds\bigg).
\]
Define $x_t$ by the (deterministic) integral equation
\begin{equation}
\label{inteq}
x_t = x_0 + \sum_{\zeta \in \NN} \zeta
\int_0^t Q_\zeta(x_s) ds
\end{equation}
and assume that it is unique for all $t \ge 0$.
Fix a time horizon $T > 0$ and let
\begin{align*}
B &:= \max_{t \le T} |x_t|,
\\
M_\zeta &:= \max_{|x|\le B} |Q_\zeta(x)|
\\
L_\zeta &:=\sup_{\substack {|x|, |y| \le B \\ x \neq y}  }
\frac{|Q_\zeta(x) - Q_\zeta(y)|}{|x-y|}
\\
\tau_\nn & :=\inf\{t>0:~ |X_{\nn,t}|>N B\}.
\end{align*}
We then have:
\begin{multline*}
\Delta_{\nn,t}:=
\frac{X_{\nn,t}}{N}-x_t =
\frac{X_{\nn,0}}{N}-x_0
+ \sum_{\zeta \in \NN} \zeta
\bigg[\frac{1}{N} \Phi_\zeta\bigg( N \int_0^t Q_\zeta(X_{\nn,s}/N) ds \bigg)
- \int_0^t Q_\zeta(X_{\nn,s}/N)ds \bigg]
\\
+ \sum_{\zeta \in \NN} \zeta \int_0^t
\big( Q_\zeta(X_{\nn,s}/N) - Q_\zeta(x_s) \big) ds
\end{multline*}
Suppose that $t \le T \wedge \tau_\nn$. Then,
for all $s \le t$,
\[
|Q_\zeta(X_{\nn,s}/N) - Q_\zeta(x_s)| \le
L_\zeta  |\Delta_{\nn,s}|.
\]
So, if we let
\[
\EE_{\nn,t}:= \frac{X_{\nn,0}}{N}-x_0
+ \sum_{\zeta \in \NN} \zeta
\frac{1}{N} 
\bigg[
\Phi_\zeta\bigg( N \int_0^t Q_\zeta(X_{\nn,s}/N) ds \bigg)
- N \int_0^t Q_\zeta(X_{\nn,s}/N)ds 
\bigg],
\]
we have, by the Gronwall-Bellman lemma, that 
\[
|\Delta_{\nn,t}| \le |\EE_{\nn,t}| 
\exp \big( t \sum_{\zeta\in \NN} |\zeta| L_\zeta\big),
\quad \text{if }
t \le T\wedge \tau_\nn.
\]
Let
\[
\Phi^*_\zeta(t) := \sup_{s \le t} |\Phi_\zeta(s)-s|.
\]
We recall that, as $N \to \infty$,
\begin{equation}
\label{plln}
\frac{1}{N} \Phi^*_\zeta(Nt) \to 0 ~ \text{ a.s.}
\end{equation}
If $s \le t \le \tau_\nn$, we have $X_{\nn,s}/N \le B$
(definition of $\tau_\nn$) and so $Q_\zeta(X_{\nn,s}/N) \le M_\zeta$,
implying that
\[
\sup_{t \le T \wedge \tau_\nn}
|\EE_{\nn,t}| \le \big|\frac{X_{\nn,0}}{N}-x_0\big| + 
\sum_{\zeta \in \NN} |\zeta| \frac{1}{N} \Phi^*_\zeta(N M_\zeta T)
\]
which converges to zero, a.s., due to \eqref{plln}.
Since
\[
\sup_{t \le T \wedge \tau_\nn}
|\Delta_{\nn,t}|
\le 
\sup_{t \le T \wedge \tau_\nn}
|\EE_{\nn,t}| 
\exp \big( T \sum_{\zeta\in \NN} |\zeta| L_\zeta\big),
\]
we have
\[
\sup_{t \le T \wedge \tau_\nn}
|\Delta_{\nn,t}| \to 0, ~ \text{ a.s.}
\]
Now observe that
\begin{align*}
P(\tau_\nn \le T) & \le P(\sup_{t \le T \wedge \tau_\nn} 
|X_{\nn,t}| > NB)
\\
& \le P(\sup_{t \le T \wedge \tau_\nn} |\Delta_{\nn,t}|
+ \sup_{t \le T \wedge \tau_\nn} |x_t| > B) \to 0.
\end{align*}
So we have $\sup_{t \le T} |\Delta_{\nn,t}|  \to 0$ a.s.

To show that $x_t$, defined via \eqref{inteq},
satisfies the ODE $\dot x = v(x)$ with
$v$ given by \eqref{vcomp} is a
matter of straightforward (but tedious) algebra,
see the Appendix.

Uniqueness and analyticity
of the solution of the ODE is immediate from the form of the vector
field (its components are polynomials of degree 2 and hence locally 
Lipschitzian).

To show that the trajectories do not explode, we 
consider the function 
\[
V(x) := \sum_A x^A.
\]
It is a matter of algebra to check that
\[
\langle \nabla V (x), v(x) \rangle 
= \sum_A v^A(x) =
\sum_A \alpha^A - \delta x^F
\]
which (since $\delta > 0$)
is negative and bounded away from zero for $x$ outside a bounded 
set of $\R_+^{\PP(F)}$  containing the origin.
We then apply the Lyapunov criterion for ODEs to conclude
that $x_t$ is defined for all $t \ge 0$ and this justifies the
fact that we could choose an arbitrary time horizon $T$ earlier in
the proof. \qed

\paragraph{\em Comment:}
The quantities defined in \eqref{phipsi},
have physical meanings as follows: 
\begin{align*}
\phi_d^A(x) &:= \sum_{B \supset A} x^B =
\text{no.\ of peers from which an $A$-peer can download},
\\
\phi_s^A(x) &:= \sum_{B \not \sim A} x^B =
\text{no.\ of peers an $A$-peer can swap with},
\\
\psi_d^A(x) &:= \sum_{a \in A} x^{A-a} =
\text{no.\ of peers which can assume label $A$ after a download},
\\
\psi_s^{A,B}(x) &:= \sum_{a \in A \cap B} x^{A-a} = 
\text{no.\ of peers which can assume label $A$ after a $B$-peer swap}.
\end{align*}
It is helpful to keep these in mind because they aid in
writing down the various parts of $v(x)$, again, see the Appendix.

\section{Behaviour of the limiting ODE}\label{ODEbeh}

Concerning the ODE $\dot x = v(x)$ we consider again three cases,
just as in the stochastic model: an open system ($\alpha^A >0$ for 
at least one $A$ and $\delta > 0$), a closed dissipative system
($\alpha^A=0$ for all $A$ and $\delta > 0$), and a closed
conservative system ($\alpha^A=0$ for all $A$ and $\delta=0$).
Qualitatively, the behaviour is different in each case.
In this section, we will try to exemplify this behaviour by means
of examples.

%
%

\subsection{The ODE in the absence of BitTorrent incentives}
Absence of BitTorrent incentives means that the file is not
split into chunks, i.e.\ $n=1$. 
Thus, the Markov chain is $X_t = (X^\varnothing_t, X^1_t :=X^F_t)$, i.e.\
two-valued $A \in\{\emptyset, 1\}$.
The rates for this case were reported earlier in \eqref{kmc}.
To find the fluid limit, we use \eqref{phipsi} and \eqref{vcomp}, 
keeping in mind the interpretation of each of the terms in \eqref{phipsi}.

For $A=\emptyset$, we have:
the number of peers from which an $\emptyset$-peer
can download from is $\phi_d^\emptyset(x) = x^1$;
the number of peers that can swap with an $\emptyset$-peer is
$\phi_s^\emptyset(x)=0$;
since no peer can assume value $\emptyset$ after an interaction,
we have $\psi_d^\emptyset(x) = \psi_s^{\emptyset,B}(x) = 0$.
Hence, in the formula for \eqref{vcomp} for $v^\emptyset(x)$
only the first parenthesis survives:
\[
v^\emptyset(x) = \alpha^\varnothing -\beta x^\emptyset x^1.
\]
For $A=1$ we have:
$\phi_d^1(x) = \phi_s^1(x)=0$, $\psi_d^1(x) = x^\emptyset$,
$\psi_s^{1,\emptyset}(x)=0$.
Thus, 
\[
v^1(x) = \alpha^1 + \beta x^\emptyset x^1 -\delta x^1.
\]
Here, by simply letting $x:=x^\varnothing$ and $y:=x^1$,
the ODE is
\begin{align*}
\dot x &= \alpha^\varnothing -\beta x y \\
\dot y &= \alpha^1+ \beta xy -\delta y.
\end{align*}

If $\alpha^\varnothing=\alpha^1=\delta=0$ (closed conservative system),
we have $x+y$=constant, say $=1$, and
\[
\dot x = -\beta x(1-x),
\]
the solution of which is the logistic function,
\begin{center}
\begin{tabular}{lcr}
$\displaystyle x_t = \frac{x_0}{x_0+(1-x_0)e^{\beta t}}$:
&
&
\begin{minipage}{7cm}
\epsfig{file=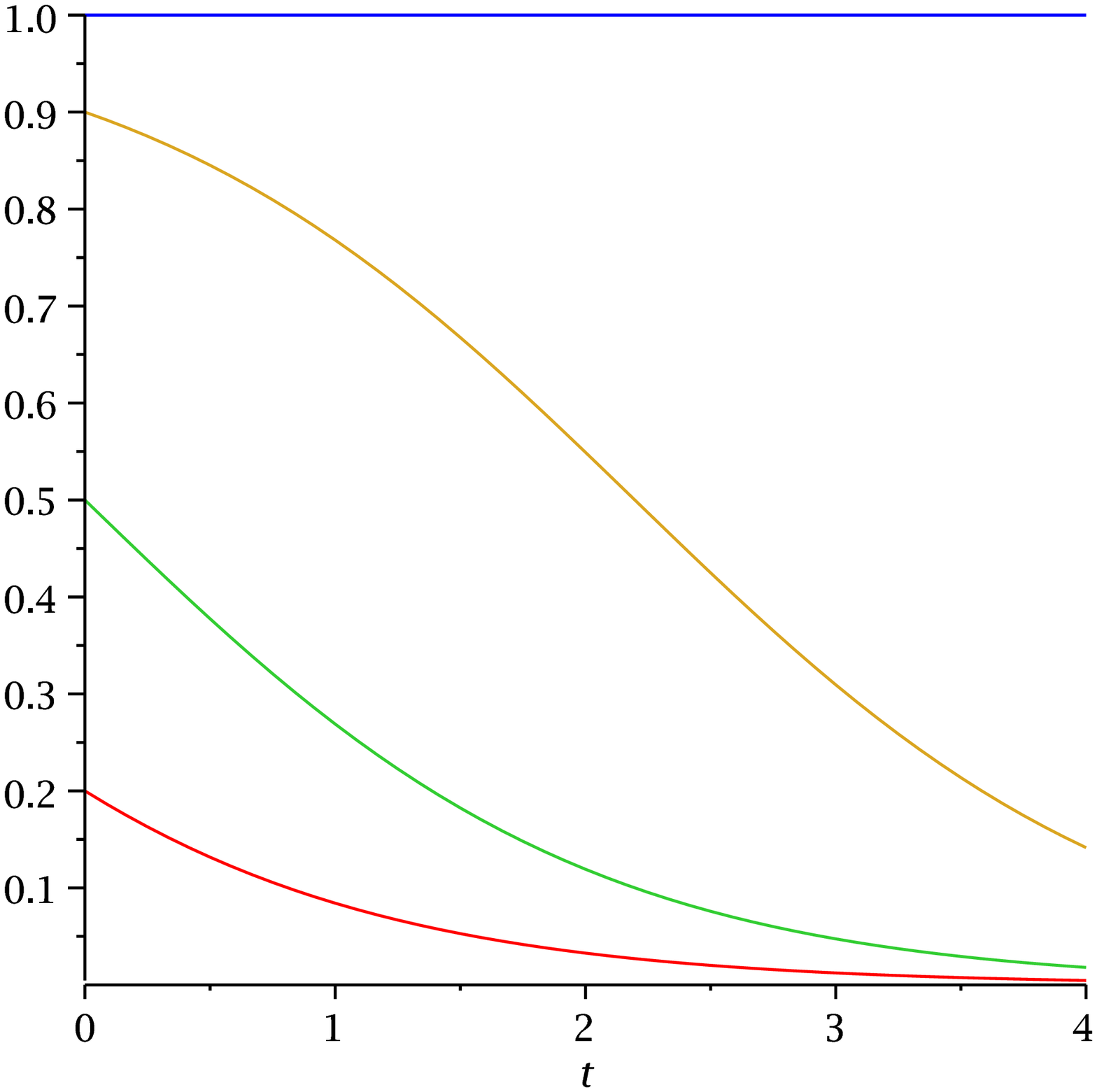,height=3.5cm}
\end{minipage}
\end{tabular}
\end{center}

If $\alpha^\varnothing=\alpha^1=0$, but $\delta >0$ (closed
dissipative system) we have the classical deterministic
SIR epidemic
\begin{align*}
\dot x &=  -\beta x y \\
\dot y &= \beta xy -\delta y
\end{align*}
the integral curves of which can be found by solving
\[
\frac{dy}{dx} = \frac{\delta}{\beta x} -1,
\]
whence it follows
\begin{center}
\begin{tabular}{lcr}
$\displaystyle 
y = (x_0+y_0) + \frac{\delta}{\beta} \log(x/x_0)-x$:
&
&
\begin{minipage}{7cm}
\epsfig{file=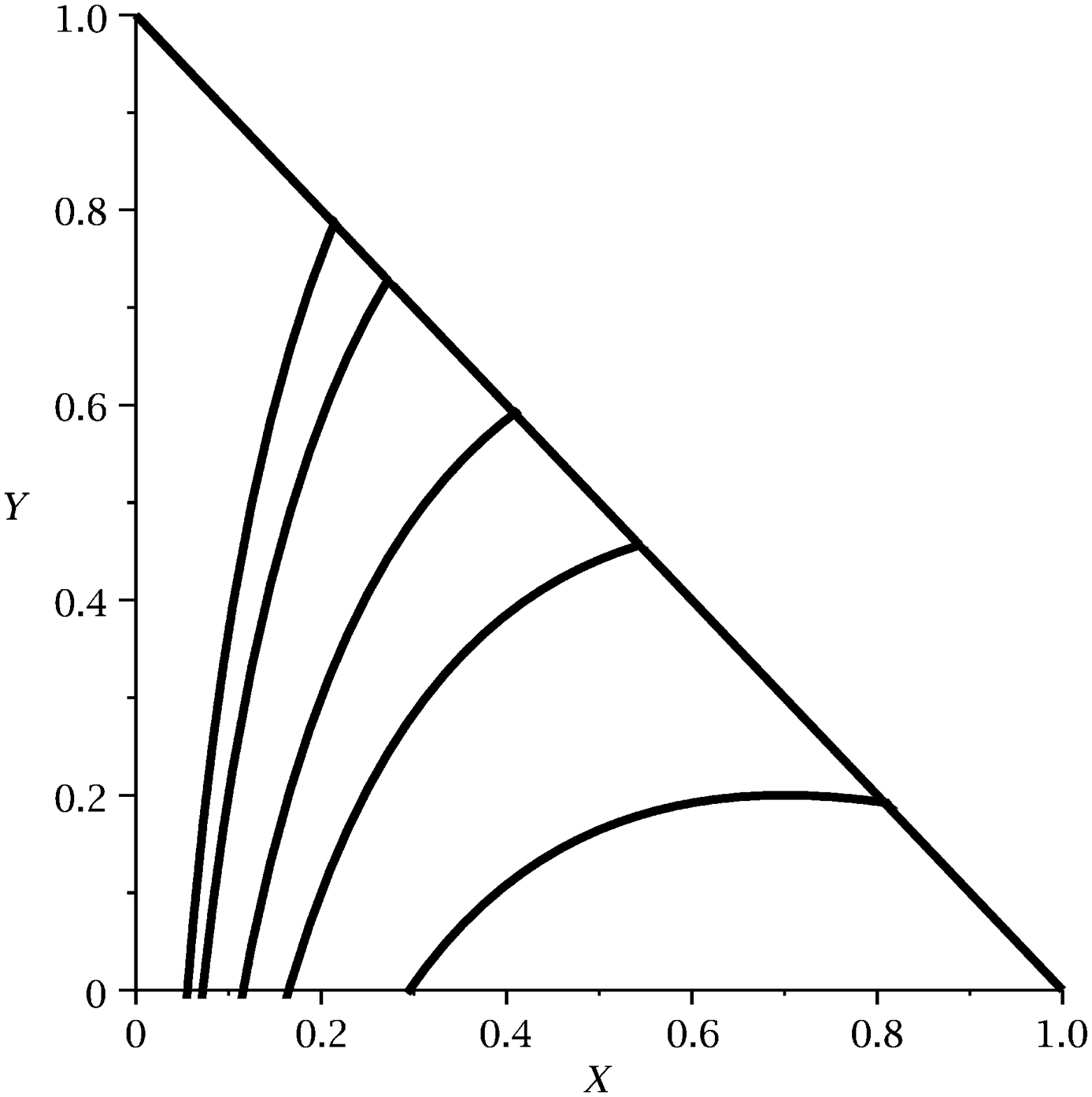,height=5cm}
\end{minipage}
\end{tabular}
\end{center}
Assume $x_0 + y_0 \le 1$.
Notice that the integral curves
are monotonic as long as they start from a point $(x_0, y_0)$
with $x_0 > \delta/\beta$.
They eventually converge to a point of the form $(x^*,0)$
where
\[
(x_0+y_0) + \frac{\delta}{\beta} \log(x^*/x_0)-x^* =0.
\]
In other words, $(x_t, y_t) \to (x^*, 0)$ as $t \to \infty$ and, in fact,
does so exponentially fast.

Notice that if the initial state is in the interior of
the positive orthant then the boundary cannot be reached in finite time.
This is in contrast to the corresponding stochastic model
which can reach the boundary in finite time with positive probability.
In fact, in the open case, it will reach the boundary in finite
time with probability $1$ due to (positive) recurrence.
This remark is generic: it applies to any dimension.

Next, consider the open system case assuming, for simplicity,
that $\alpha^\varnothing = \lambda >0$, but $\alpha^1=0$.
We have
\begin{align*}
\dot x &= \lambda -\beta x y \\
\dot y &= \beta xy -\delta y.
\end{align*}
Since $\delta >0$, the Lyapunov function argument shows that
the system is asymptotically stable. In fact, there is
a unique asymptotic equilibrium which attracts all initial conditions
$(x_0, y_0)$ with $y_0 > 0$.
This unique equilibrium is given by
\[
x^* = \frac{\delta}{\beta}, ~ y^* = \frac{\lambda}{\delta},
\]
as can be seen by setting the right hand side of the ODE equal to zero
\cite{Daley-Gani}.
It should be noticed that the trajectories can be spirals around 
$(x^*, y^*)$.
A typical situation 
is shown below.
\begin{center}
\epsfig{file=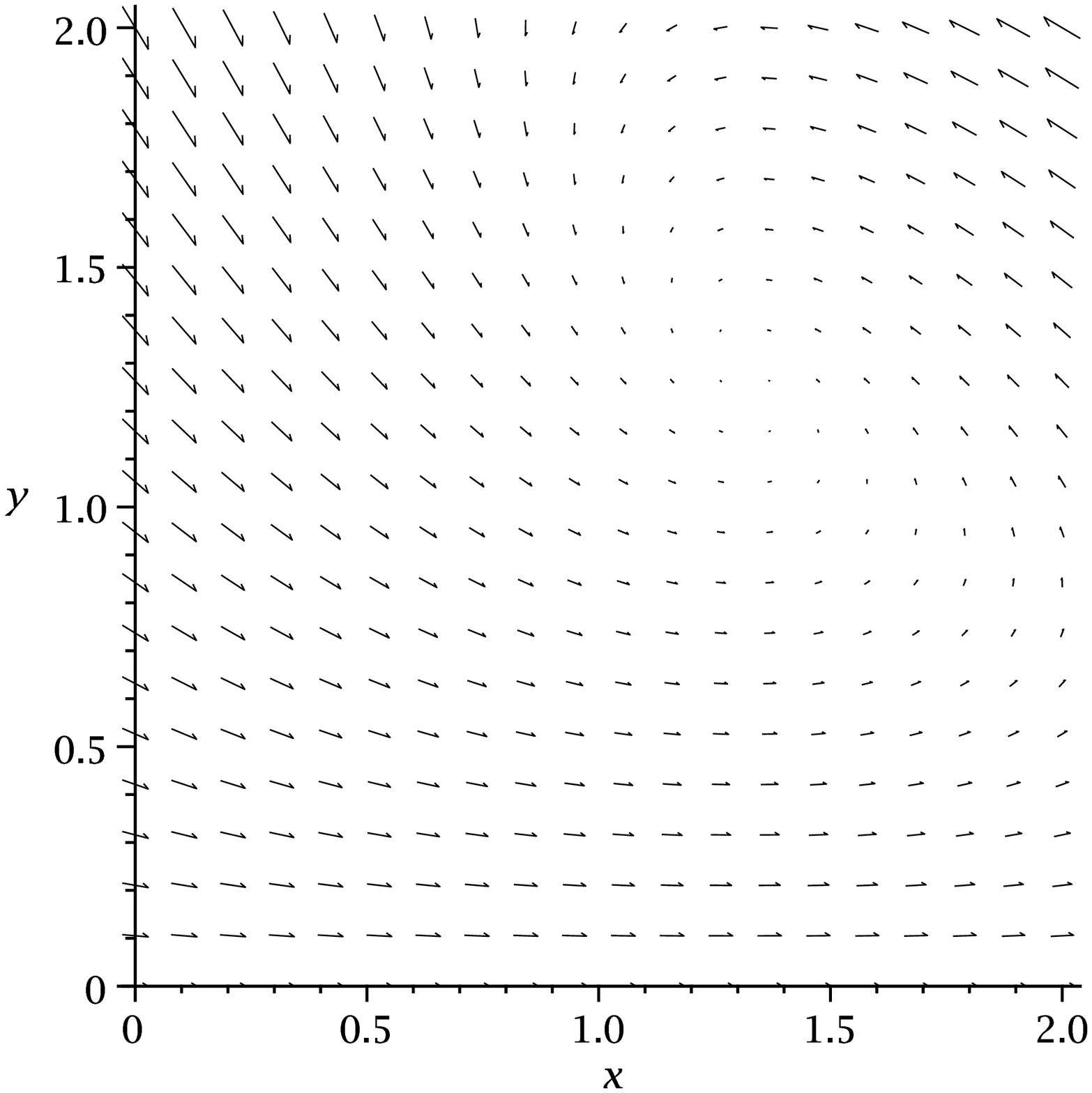, height=6cm}
\end{center}
In this vector field plot, we took $\beta=3$, $\lambda=5$, $\delta=4$,
so $(x^*, y^*) = (1.33, 1.25)$.
If $y_0=0$ then the trajectory converges to infinity. Indeed,
if no peers initially possess the file and no such peers ever show
up, then the only thing that can happen is an accumulation, at rate 
$\lambda$, of ever demanding peers.
The remedy is, clearly, the imposition of $\alpha^1>0$ no matter
how small; then $(x^*, y^*)$ is globally attracting for
all initial states in the closed positive orthant.

In this last example, we find that the eigenvalues of the differential
of the vector field at $(x^*, y^*)$ are complex conjugates
if and only if $\lambda\beta < 4 \delta^2$, and so this is the condition
for the spiralling of trajectories.

\subsection{The ODE for $n=2$ chunks}

Here $A$ can take $4$ values: $\emptyset$, $\{1\}$, $\{2\}$, $\{1,2\}:= F$.
We work out the expression for each
component separately.
Since an $\emptyset$-peer can only download from
$\phi_d^\emptyset(x)=x^1+x^2+x^{12}$
peers\footnote{Where $x^{12}=x^F$ and indexes 1 and 2 represent the chunks that
form a partition of $F$.}, but no peer can download from it or swap with it,
we have $\phi_s^\emptyset(x)=\psi_d^\emptyset(x)=\psi_s^\emptyset(x)=0$, and
so 
\[
v^\emptyset(x) = \alpha^\varnothing -\beta x^\emptyset (x^1+x^2+x^{12}).
\]
A $1$-peer can download from $\phi_d^1(x)=x^{12}$ peers;
it can swap with $\phi_s^1(x)=x^2$ peers;
the number of peers that take value $1$ after a download is
$\psi_d^1(x)+x^\emptyset$; and, since no peer can take the
label $1$ after a swap, we have $\psi_s^{1,B}(x)=0$ $\forall B$.
We thus have
\begin{align*}
v^1(x) &= \alpha^1 -x^1(\beta x^{12} + \gamma x^2)
+ \beta \psi_d^1(x) \sum_{B: 1 \in B} \frac{x^B}{1+|B\setminus \{1\}|}
\\
&= 
\alpha^1 -x^1(\beta x^{12} + \gamma x^2)
+ \beta  x^\emptyset(x^1 + \tfrac{1}{2} x^{12}).
\end{align*}
The expression for $v^2(x)$ is symmetric to that of $v^1(x)$.
A $12$-peer can neither download from or swap with anyone,
so $\phi_d^{12}(x) = \phi_s^{12}(x)=0$;
there are $\psi_d^{12}(x) = x^1+x^2$ peers which can take label $12$
after a download;
since only a $2$-peer can assume label $12$ after a swap, we have
$\psi_s^{12, 1}(x) = x^2$; and, for the same reason,
$\psi_s^{12, 2}(x) = x^1$.
Thus,
\begin{align*}
v^{12}(x) &= \alpha^{12} + \beta \psi_d^{12}(x) x^{12} + 
\frac{\gamma \psi_s^{12,1}(x) x^1}{1+|\{1\}\setminus \{1,2\}|} 
+  
\frac{\gamma \psi_s^{12,2}(x) x^2}{1+|\{2\}\setminus \{1,2\}|} 
-\delta x^{12}
\\
&=
\alpha^{12} + \beta (x^1+x^2) x^{12} 
+ \gamma (x^2 x^1 +  x^1 x^2)
-\delta x^{12}.
\end{align*}
We thus have,
\begin{align*}
&\dot x^\emptyset 
=  
\alpha^\varnothing -\beta x^\emptyset (x^1+x^2+x^{12}) 
\\
&\dot x^1 
=  
\alpha^1 -x^1 (\beta x^{12} +\gamma x^2)
+ \beta x^\emptyset ( x^1 + \tfrac{1}{2}x^{12} ) 
\\
&\dot x^2
=
\alpha^2 -x^2 (\beta x^{12} +\gamma x^1)
+ \beta x^\emptyset ( x^2 + \tfrac{1}{2}x^{12} ) 
\\
&\dot x^{12}
=
\alpha^{12} + \beta (x^1+x^2) x^{12} + 2\gamma x^1 x^2 -\delta x^{12}.
\end{align*}

\paragraph{Case 1: closed conservative system.}
Consider the $n=2$ case again, and assume that $\alpha^A =0$ for all $A$,
$\delta=0$ (a closed conservative system). Assume further that
$\gamma=0$. Let 
\[
x=x^\varnothing,\quad u=x^1+x^2,\quad w=x^{12}.
\]
We see a reduction in dimension from $4$ to $3$ (owing to
that fact that the sum of the coordinates is constant)
and a further reduction from $3$ to $2$ (owing to the
fact that, for $\gamma=0$, the vector field depends
on $x^1$, $x^2$ through their sum).
We have:
\begin{align*}
\dot x &= -\beta x(u+w) \\
\dot u &= -\beta uw + \beta x(u+w) \\
\dot w &= \beta uw.
\end{align*}
On assuming that 
\[
x+u+w=1
\]
and eliminating the variable $w$ we obtain
\begin{align*}
\dot x &= -\beta x(1-x) \\
\dot u &= \beta u^2 - \beta u(1-x) + \beta x(1-x).
\end{align*}
We can qualitatively see the behaviour of this ODE by
looking at the vector field in the $x-u$ plane. A typical
picture is as follows:
\begin{center}
\epsfig{file=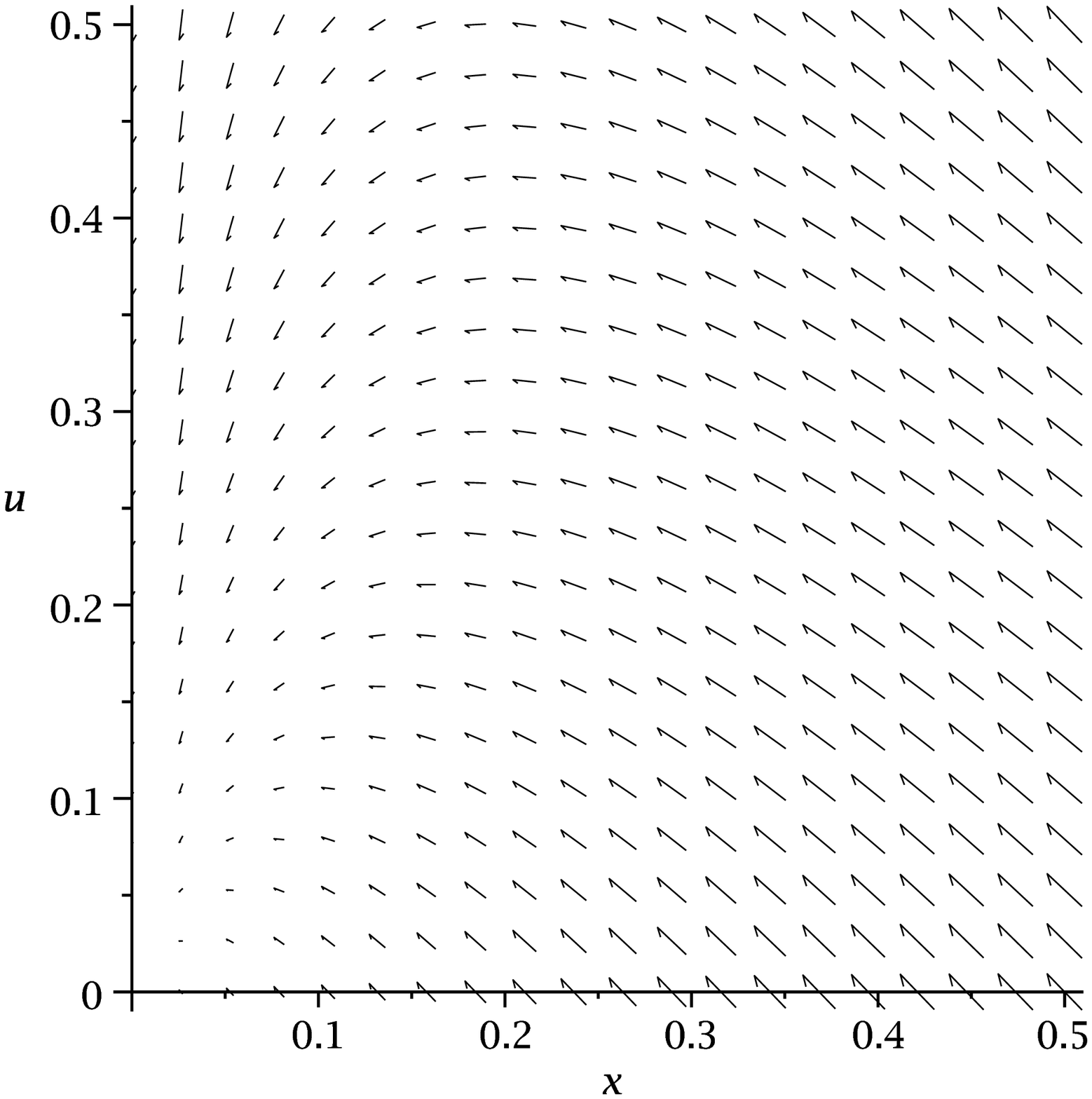, height=6cm}
\end{center}
The first is an autonomous equation, encountered earlier.
Its solution is $x_t = x_0/(x_0+(1-x_0)^{\beta t})$.
Letting
\[
f := 1-w = x+u
\]
we have
\[
\dot f = \beta f (1-f) + \beta (1-f) x_t
\]
Assume that $w_0 > 0$, i.e.\ $f_0 < 1$.
The solution of the last ODE is:
\[
f_t = \frac{f_0(1-f_0)^{-1} + x_0 \beta t}
{x_0+(1-x_0)e^{\beta t} + f_0(1-f_0)^{-1} + x_0 \beta t}.
\]
Hence
\[
w_t = \frac{x_0 + (1-x_0) e^{\beta t}}
{x_0 + (1-x_0) e^{\beta t} + x_0 \beta t + (1-w_0)w_0^{-1}}.
\]
We have $\lim_{t \to \infty} w_t =1$,
as expected.
{From} this, we can estimate the time $\tau=\tau(x_0,w_0, \epsilon)$
required for the 
system to reach the set $\{(x,u,w):~ w > 1-\epsilon\}$, 
where $0< \epsilon <1$ is a given (small) number.
This time can be translated into an estimate for the expected time
required for the original stochastic system to be absorbed by
the state where all peers possess the full file. We shall not attempt
the justification of this statement here.
The time $\tau$ is the solution of the transcendental equation
\[
\frac{1-w_0}{w_0} + x_0 \beta \tau
= \frac{\epsilon}{1-\epsilon}(x_0+(1-x_0) e^{\beta \tau}).
\]
Typical behaviour of this time as a function of
$x_0$ is as follows:
\begin{center}
\epsfig{file=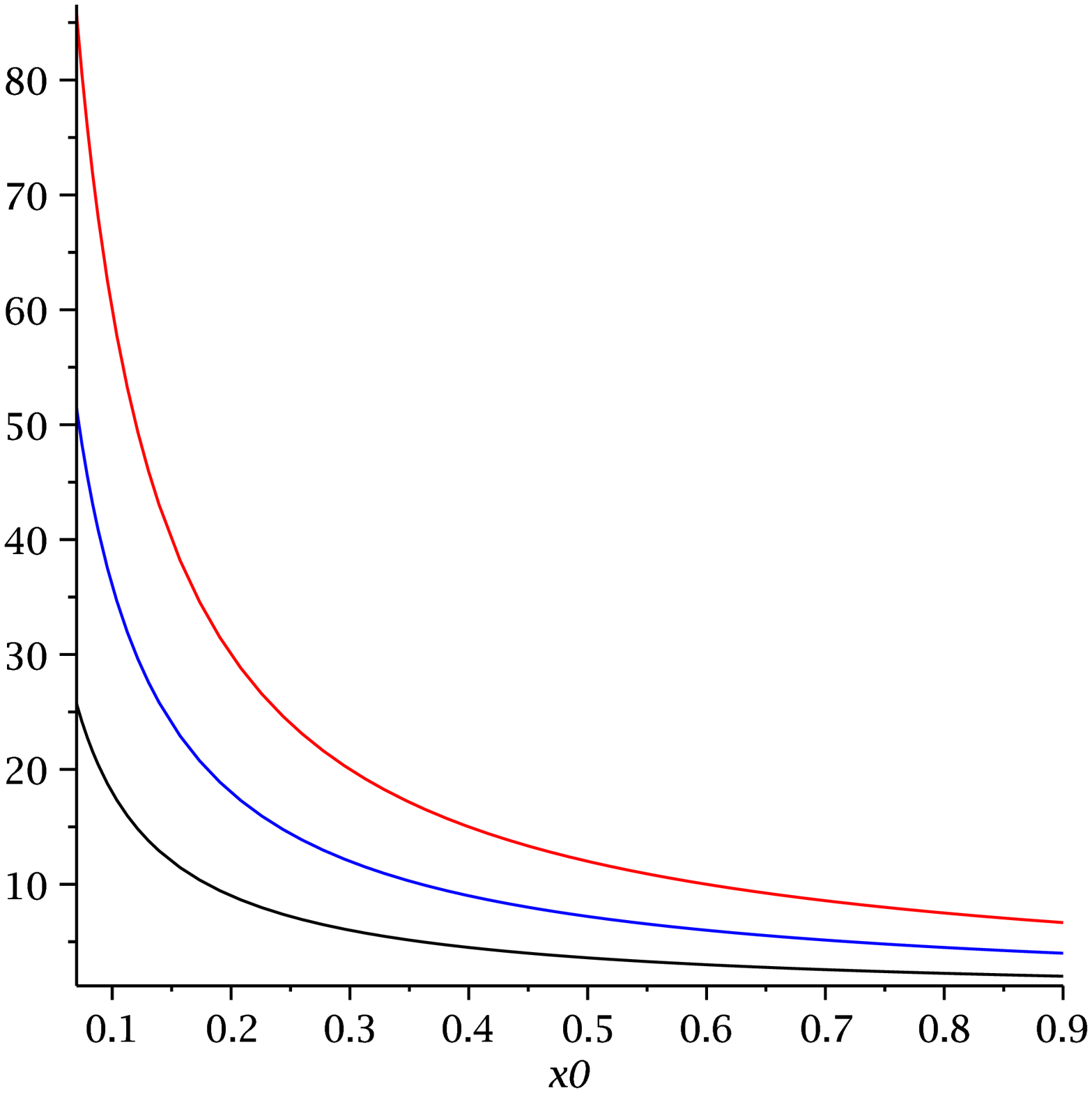, height=6cm}
\end{center}
The three graphs correspond to varying values of $\beta$
ranging from $1$ (top curve) to $5$ (bottom curve) and to 
$\epsilon=0.001$, $w_0=0.1$.

\paragraph{Case 2: closed dissipative system.}
Consider again the $n=2$ case with no arrivals, with swap rate $\gamma$ equal
to zero, but with departure rate $\delta > 0$. Using the same
notation as above, we have
\begin{align*}
\dot x &= -\beta x(u+w) \\
\dot u &= -\beta uw + \beta x(u+w) \\
\dot w &= \beta uw-\delta w.
\end{align*}
To substitute out one parameter, change the time 
variable to $s=\beta t$, let
$\rho=\delta/\beta$, and write $x'$ for $dx/ds$, so that
\begin{align*}
x' &= -x(u+w) \\
u' &= -uw + x(u+w) \\
w' &= uw-\rho w.
\end{align*}

Assume
\[
x_0+u_0+w_0=1,
\]
so $x_t+u_t+w_t < 1$ for all $t > 0$.
We cannot eliminate the variable $w$ now since there is no
obvious conserved quantity, but we can study the
equilibria of the system.
It seems that the ODE has no obvious analytical solution.
However it can either be integrated numerically or its solution in
terms of power series can be easily found.
Setting the right hand side equal to zero we see that, necessarily,
$w=0$, which leaves us with the possibility $xu=0$.
So points of the form
\[
(x,u,w) = (x,0,0), ~ (0,u,0)
\]
are equilibria.
But not all of them are stable.
For example, any point of the form $(0,u,0)$ with $u>\rho$
is an unstable equilibrium.
Specifically, the differential of the vector field
is given by
\[
Dv =
\begin{pmatrix}
-(u+w) &  -x   & -x  
\\
u+w   &  x-w   & x-u
\\
0    & w     & u-\rho
\end{pmatrix}
\]
Evaluated at $(0,u,0)$, it gives
$Dv =
\textstyle
\begin{pmatrix}
0 & 0 & 0
\\
0 & 0 & -u
\\
0    & 0     & u-\rho
\end{pmatrix}$,
and this has eigenvalues $0, -u, u-\rho$.
Evaluated at $(x,0,0)$,
it gives
$Dv=
\textstyle
\begin{pmatrix}
0 & -x & -x
\\
0 & x & x
\\
0    & 0     & \rho
\end{pmatrix}$,
and this has eigenvalues $0, x, -\rho$.
Therefore, the only stable equilibria
are points of the form
\[
(x,u,w) = (0,u,0), \quad 0 \le u < \rho.
\]
In terms of the original variables, the stable equilibria are
\[
(x^\varnothing, x^1, x^2, x^{12}) = (0, x^1, x^2, 0),
\quad 0 \le x^1+x^2 < \rho.
\]
This is as expected: since there is no swapping ($\gamma=0$),
the system eventually settles to a situation where there
are peers with label $1$ and peers with label $2$.
Had $\gamma$ been positive, $x^1$, $x^2$ could not have simultaneously
been positive in equilibrium.

\paragraph{Case 3: Open system.}
Consider the situation as in Case 2, but add arrivals of peers 
(known as seeds) possessing the full file.
Choosing variables appropriately, we have
\begin{align*}
x' &= -x(u+w) \\
u' &= -uw + x(u+w) \\
w' &= \lambda+ uw-\rho w.
\end{align*}
In terms of the original variables, we here have $\alpha^\varnothing=
\alpha^1=\alpha^2=0$, $\alpha^{12} = \lambda \beta > 0$.
Here we see that the system eventually settles to the stable
equilibrium
\[
(x,u,w) = (0,0, \lambda/\rho).
\]
We can easily see that the eigenvalues of the differential
of the vector field at the stable equilibrium
are both real and negative: $-\lambda/\rho$ and $-\rho$
(the first one has algebraic multiplicity 1 but geometric multiplicity 2).
Hence there is no possibility of spiralling.
\begin{center}
\epsfig{file=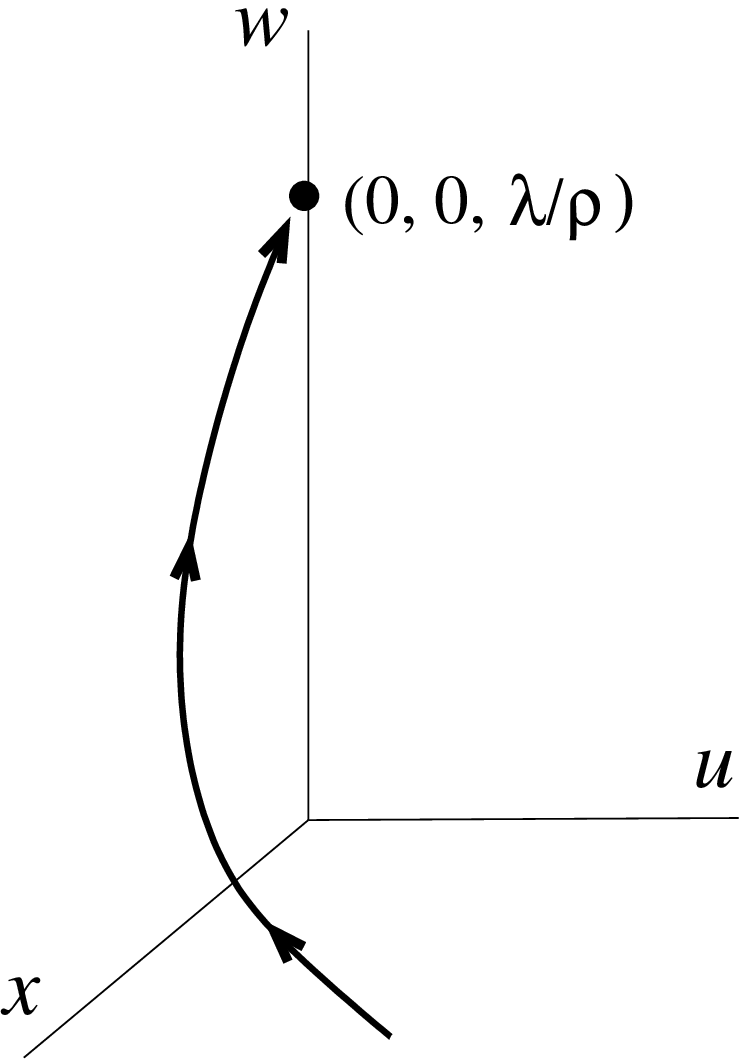, height=5.5cm}
\end{center}

\subsection{Time to settle}
\label{settle}
Estimating the time for the system to reach an equilibrium
requires cooking up an appropriate Lyapunov function.
The obvious Lyapunov function used earlier gives a crude
lower bound.

Consider a general deterministic system (open or closed), 
i.e.\ the differential equation $\dot x = v(x)$ with $v(x)$
given by \eqref{vcomp}. Let $V(x) = |x| = \sum_{A \subset F} x^A$.
Then
\[
\dot V = |\alpha| -\delta x^F,
\]
where $|\alpha| = \sum_{A \subset F} \alpha^A$ is the total arrival rate.
Using the naive inequality $x^F \le V$ we obtain the following:

\begin{corollary}
If the system starts from $x_0$ and if $x_t \to x^*$ then
the time $\tau_r$ required for the trajectory to reach an $r$-ball
centred at $x^*$ satisfies
\[
\tau_r \ge \frac{1}{\delta} \log \frac{|x_0|+|\alpha|}{|x^*|+r}.
\]
\end{corollary}

In the other direction, 
consider the last component $v^F(x)$ of the vector field.
Since no seed ($F$-peer) can download or swap, we have
$\phi_d^F(x)=\phi_s^F(x)=0$. So 
\begin{align*}
v^F(x) &= \alpha^F 
+ \beta \psi_d^F(x) x^F
+ \gamma \sum_{B \subset F} \psi_s^{F,B}(x) x^B - \delta x^F
\\
&= \alpha^F + \beta \sum_{j \in F} x^F x^{F-j}
+ \gamma \sum_{B \subset F} \sum_{i \in B} x^B x^{F-i}
 - \delta x^F
\\
&=: \alpha^F + v_+^F(x) - \delta x^F,
\end{align*}
where $v_+^F(x)$ is a quadratic form with positive coefficients.

Assume that the system is closed, so that $\alpha^F=0$ in particular.
Let $\overline {v_+^F}$
be an upper bound on $v_+^F(x)$, i.e., a (good) bound that depends 
on the initial state
$x_0$ and the parameters $\beta, \gamma$. 

\begin{corollary}
If  the system is closed and
$\overline {v_+^F} < \delta$,
then 
\[
\tau_r \le \frac{1}{\overline {v_+^F}-\delta} \log\frac{x_0^F}{r}.
\]
\end{corollary}

We conjecture that $v_+^F(x_t)$ decreases along the trajectory $x_t$
as long as $v_+^F(x_0) < \delta$.
If so, the bound above is 
$\tfrac{1}{v_+^F(x_0)-\delta} \log\frac{x_0^F}{r}$.

\section{An example of the evaluation of performance improvement in presence of
BitTorrent incentives}\label{incentives-sec}

We address the following question: When is it advantageous to split a file 
into chunks? In other words, assuming we fix certain system parameters (e.g.,
arrival rates), will peers acquire the file
faster if the file is split into chunks?
We attempt here to answer the question in a simple case only
by using the deterministic approximation.
Let $\lambda$ be the total peer arrival rate. Let $\beta$ be the download rate.
Assume that only $\varnothing$ peers arrive exogenously.
In the absence of BitTorrent incentives, we have the single-chunk case
\begin{align*}
\dot x^\varnothing &= \lambda - \beta x^\varnothing x^1
\\
\dot x^1 &= \beta x^\varnothing x^1 -\delta x^1.
\end{align*}
The globally attracting stable equilibrium is given by 
\[
x^* = (\delta/\beta, ~ \lambda/\delta).
\]
Consider splitting into $n=2$ chunks.
Let $\widetilde x$ be the state of the system. 
Suppose that the new parameters are 
$\ttt \lambda = \lambda$, $\ttt \delta = \delta$, $\ttt \beta$, $\ttt \gamma$.
Then
\begin{align*}
&\dot  {\ttt x}^\emptyset
=
\lambda -\ttt \beta \ttt x^\emptyset (\ttt x^1+\ttt x^2+\ttt x^{12})
\\
&\dot {\ttt x}^1
=
-\ttt x^1 (\ttt \beta \ttt x^{12} +\ttt \gamma \ttt x^2)
+ \ttt \beta \ttt x^\emptyset ( \ttt x^1 + \tfrac{1}{2}\ttt x^{12} )
\\
&\dot {\ttt x}^2
=
-\ttt x^2 (\ttt \beta \ttt x^{12} +\ttt \gamma \ttt x^1)
+ \ttt \beta \ttt x^\emptyset ( \ttt x^2 + \tfrac{1}{2}\ttt x^{12} )
\\
&\dot {\ttt x}^{12}
=
\ttt \beta (\ttt x^1+\ttt x^2) \ttt x^{12} + 2\ttt \gamma \ttt x^1 \ttt x^2 -\delta \ttt x^{12}.
\end{align*}
The new equilibrium is easily found to be
\[
\ttt x^* = \bigg(\frac{\delta}{\ttt \beta} \big(\frac{\delta}{\lambda}u+1 \big)^{-1},~
\frac{u}{2},~ \frac{u}{2},~ \frac{\lambda}{\delta}\bigg),
\]
where $u$ is the positive number which solves
\[
q(u)=0,
\]
and
\begin{equation}
\label{qqq}
q(u):= u^2 + \frac{2 \ttt \beta\lambda}{\ttt \gamma \delta} u 
- \frac{2 \lambda}{\ttt \gamma}.
\end{equation}
To see this, set the vector field equal to zero and solve for
$\ttt x^A$, $A\in\{\varnothing, 1, 2, 12:= F\}$.
The simplest way to obtain the solution is by first adding the equations;
this gives
\[
\lambda-\delta \ttt x^{12} = 0,
\]
whence $\ttt x^{12} = \lambda/\delta$.
Then add the middle two equations after setting $u=\ttt x^1 + \ttt x^2$:
\[
-\ttt \beta u \ttt x^{12} -2\ttt \gamma \ttt x^1 \ttt x^2 + 
\ttt \beta x^\varnothing 
(u+ \ttt x^{12}) =0.
\]
Replace $\ttt x^{12}$ by $\lambda/\delta$ and observe that, due to symmetry,
$\ttt x^1 = \ttt x^2 = u/2$. This gives the quadratic equation $q(u)=0$
with $q$ defined by \eqref{qqq}. 
Finally, the first equation becomes 
\[
\lambda - \ttt \beta \ttt x^\varnothing (u + \lambda/\delta)=0,
\]
which is solved for $\ttt x^\varnothing$ giving:
\[
\ttt x^{*\varnothing} = \frac{\lambda}{\ttt \beta}
\big(u+\frac{\lambda}{\delta})^{-1}
= \frac{\delta}{\ttt \beta} \big(\frac{\delta}{\lambda}u+1 \big)^{-1}
< \frac{\delta}{\ttt \beta}.
\]
Thus:

\begin{corollary}
If $\ttt \beta \ge \beta$ then
\[
{\ttt x}^{*\varnothing}< x^{*\varnothing}  .
\]
\end{corollary}

So, by introducing splitting into chunks, we have fewer peers
who have no parts of the file at all.
Using Little's theorem (see below), this can be translated into smaller waiting
time from the time a peer arrives until he gets his first chunk.

Suppose now we are interested in determining how long it
will take for a newly arrived peer to acquire the full file.
On the average, a peer spends time equal to $\lambda^{-1} |x^*|$ 
before it exits the system. During last part of his sojourn interval
(which is a random variable with mean $1/\delta$), the peer
possess the full file. It thus takes on the average $\lambda^{-1} |x^*|-
\delta^{-1}$ for a peer to acquire the full file. Since
we assume that $\ttt \lambda = \lambda$, $\ttt \delta = \delta$, 
it suffices to show that
\[
|x^*| > |\ttt x^*|.
\]
But
\begin{align*}
|x^*| - |\ttt x^*|  
&= \left[\frac{\delta}{\beta}+\frac{\lambda}{\delta}\right]
- \left[\frac{\delta}{\ttt \beta}
\big(\frac{\delta}{\lambda}u+1 \big)^{-1} + u +\frac{\lambda}{\delta}\right]
\\
&= \frac{\delta}{\beta} 
- \frac{\delta}{\ttt \beta}\big(\frac{\delta}{\lambda}u+1 \big)^{-1} 
-u -\frac{\lambda}{\delta}
\\
&= \big(\frac{\delta}{\lambda}u+1 \big)
\left[
\big(\frac{\delta}{\beta}-\frac{\delta}{\ttt \beta}\big)
+\big(\frac{\delta^2}{\beta \lambda}-1\big) u
-\frac{\delta}{\lambda} u^2
\right]
\end{align*}
recalling $u>0$ solves $q(u)=0$.
So, $|x^*| - |\ttt x^*| > 0$ if and only if
\begin{equation}
\label{tqqq}
0> \ttt q(u) := u^2 - \big( \frac{\delta}{\beta}-\frac{\lambda}{\delta} \big) u
- \big( \frac{\lambda}{\beta} - \frac{\lambda}{\ttt \beta} \big).
\end{equation}
Define $\ttt u$ as the unique positive number which satisfies 
\[
\ttt q (\ttt u) =0.
\]

\begin{corollary}\label{incentives-coro}
If $\ttt \beta \ge \beta$,
a necessary and sufficient condition for $|x^*| > |\ttt x^*|$
is $u < \ttt u$.
\end{corollary}

This gives a set of non-vacuous conditions for achieving improvement of
performance by the introduction of BitTorrent incentives.\footnote{The
inequality conditions were mistakenly reversed in the corresponding results of
\cite{ICASSP07}.}
It can be proved that if the parameters
$\beta$, $\delta$, $\ttt \beta$, $\ttt \gamma$ are fixed
and if
$\ttt \beta \ge \beta$
then there exists a value $\lambda_0$ such that for all 
$\lambda < \lambda_0$ the inequality $u < \ttt u$ holds.
To prove this, we observe that
\[
u < \ttt u \iff q(\ttt u) > 0
\]
and study the behaviour of $q(\ttt u)$ as a function of $\lambda$ in
a neighbourhood of zero.

We conjecture that an algebraic condition involving quadratics like
$q$ and $\ttt q$ is valid for larger values of $n$ also.

To justify the use of deterministic approximation for
estimating performance measures, and, specifically, the use of mean
values, we need to show that as $N \to \infty$,
we can approximate stationary averages in the original stochastic
network by equilibria of the resulting ODE.
It is easy to show that the a.s.\ convergence to the ODE limit 
can be translated into convergence of the means, using
a uniform integrability argument. Namely,
\[
\frac{1}{N} E X_{\nn,t}
\xrightarrow[\nn \to \infty]{}  x_t \xrightarrow[t \to \infty]{} x^*,
\]
where the second limit concerns the behaviour of the ODE alone.
On the other hand, if we fix $N$ and look at the
asymptotic behaviour of the process $\tfrac{1}{N} X_{\nn,t}$ as 
$t\to \infty$, we have
\[
\frac{1}{N} E X_{\nn,t}
 \xrightarrow[t \to \infty]{}
\frac{1}{N} E \ttt X_N,
\]
where the law of $\ttt X_N$ is the stationary distribution of
the chain $(\tfrac{1}{N} X_{\nn,t})_{t \ge 0}$.
It can be proved that $\tfrac{1}{N} E \ttt X_N \to
x^*$, as $N \to \infty$.
Arguments for this will be considered in future work.
More detailed estimates on the discrepancy 
between the stochastic and deterministic systems can be found in
the recent survey paper \cite{DN08}.

We can also explain the use of $|x^*|$ as a measure of the
sojourn time in the system of a peer, by first using the
approximation outlined above and then appealing to Little's law.
This is as follows.

Consider an open \bt{X_{0}}{n}{\alpha}{\beta}{\gamma}{\delta},
i.e.\ $|\alpha| > 0$, $\delta > 0$.
We know that the Markov chain $(X_t)$ is positive recurrent and has thus
a unique stationary distribution.
It makes sense to assess the performance of the network by looking at
steady-state performance measures, such as the mean time it takes
for an $\varnothing$-peer to become an $F$-peer (a seed).
Consider then the process $(\widetilde X_t, t \in \R)$
defined to be a stationary Markov process with time index $\R$
and transition rates as those of $(X_t)$.
The law of the process $(\widetilde X_t, t \in \R)$ is unique.
Let $T^A_k, k \in \Z$ be the times at which $A$-peers arrive
(and, say, $T^A_0 \le 0 < T^A_1$, by convention). These are the 
points of a stationary Poisson process in $\R$ with rate $\alpha^A$.
Let $W^A_k$ be the sojourn time in the system of a peer
arriving at time $T^A_k$. Since, by assumption, a peer departs only
after it has acquired the full set, the time $W^A_k$ is the sum of the
times it takes for the peer to become a seed plus the time that the
peer hangs out in the system after becoming a seed (the latter is
an exponential time with mean $1/\delta$).
Clearly then, for all $t \in \R$,
\[
\sum_{B \supset A} \widetilde X^B_t = \sum_{k \in \Z}
\1(T_k^A \le t < T_k^A + W_k^A).
\]
Using Campbell's formula, we obtain
\begin{equation}
\label{ll}
\sum_{B \supset A} E \widetilde  X^B_0
= \alpha^A E^A W^A_0,
\end{equation}
where $E^A$ is expectation with respect to $P^A$--the
Palm probability of $P$ with respect to the point process $(T^A_k, k \in \Z)$.

In particular, with $A=\varnothing$, and $\lambda=\alpha^\varnothing$,
we have that
\[
E^\varnothing W_0^\varnothing = \frac{1}{\lambda}
E |\ttt X_0|,
\]
which can be read as: the mean sojourn time of a $\varnothing$-peer
is, in steady state, equal to the mean number of peers in the
system divided by the rate of arrivals of $\varnothing$-peers.
If $N$ is a parameter of the process as in Theorem \ref{ODEapprox}
then, $\lambda$ being proportional to $N$, we have that
the right hand side converges to something that is proportional
to $|x^*|$, as required.

%
%

\section{Diffusion approximation}\label{diff-sec}
Using the functional central limit theorem
for Poisson processes, we can prove, by standard methods, the following:
Again consider the sequence
\bt{X_{\nn,0}}{n} {N\alpha}{\frac{\beta}{N}}
{\frac{\gamma}{N}}{\delta} for
$N \in \N$, and let $X_{\nn,t}$ be the corresponding jump Markov chain.
Let $(x_t, t \ge 0)$ be the solution
to the ODE $\dot x = v(x)$ with initial condition $x_0$.
Let 
\[
Y_{\nn,t} := \sqrt{N} (X_{\nn,t}/N - x_t).
\]
Let $W_\zeta$, $\zeta \in \NN$, be i.i.d.\ standard Brownian motions
in $\R$. 
Finally, define the (time-inhomogeneous) Gaussian diffusion process $Y$ by
\[
dY_t = \sum_{\zeta \in \NN} \zeta \sqrt{Q_\zeta(x_t)} dW_{\zeta,t}
+ Dv(x_t) Y_t dt,
\]
where $Dv(x)$ is the matrix of partial derivatives
of $v(x)$.

\begin{theorem}
\label{diffapprox}
If $\sqrt{N} (X_{\nn,t}/N-x_0) \to 0$ as $N \to \infty$,
where $x_0 \in \R_+^{\PP(F)}$, 
then the law of $Y_\nn$ (as a sequence of probability
measures in $D[0,\infty)$ with the topology of
uniform convergence on compacta) converges weakly
to the law of $Y$.
\end{theorem}

The proof of this theorem is omitted but the reader is referred to 
\cite{Kurtz81} for the relevant arguments.

\section{Final remarks, open problems and future work}\label{concl-sec}

\subsection{Rates of convergence}
We can obtain a computable
 rate of convergence of the stochastic model to the ODE
by using a combination of large deviations techniques with 
the solution of two optimisation problems. The idea is basically
implicit in the proof of Theorem \ref{ODEapprox} and this is the
reason we wrote the proof explicitly in terms
of the driving Poisson processes $\Phi_\zeta$.

The first problem is so that we obtain
an estimate of the maximum value of $M_\zeta$ of $Q_\zeta(x)$.
In the case of a closed network, this is a quadratic optimisation
problem over the polyhedron $\{|x| \le 1\}$.

The second optimisation problem is for an estimate for $L_\zeta$,
which can be translated to an estimate for the norm of the gradient 
$\nabla Q_\zeta(x)$. In the case of a closed
network, we can estimate this by solving
a (large) number of linear programming problems.
Of course, only estimates are needed.

\subsection{Conjectures}
The first one concerns the behaviour of $v_+^F(x_t)$, and was stated
at the end of Section \ref{ODEbeh}.
The second concerning generalisation of Corollary \ref{incentives-coro}, 
was stated  in Section \ref{incentives-sec}.

The third conjecture is more vague: it basically says that we
can evaluate the performance improvement by using a large number
of chunks (say 100), by solving a number of quadratic inequalities.
To this end, it should be remarked that in a deterministic open
network, the unique equilibrium $x^*$ can be found by solving
$n+1$ equations in $n+1$ unknowns, provided that the rate of
arrivals of $A$-peers depends on $A$ through its cardinality alone:
\[
\alpha^A = \lambda_k \text{ if } |A|=k \text{ chunks.}.
\]
Indeed, by symmetry of the vector field, we see that
\[
x^{*A} = x^{*B} \text{ if } |A|=|B|.
\]
So if we define 
\[
z^k := \sum_{|A|=k} x^{*A}
\]
we will have 
\[
x^{*A} = \binom{n}{k}^{-1} z^k, \text{ if } |A|=k.
\]
Hence if we let
\[
V^k(z^0,\ldots,z^n) := \sum_{|A|=k} v^A(x^*),
\]
where $z$ and $x^*$ are related as above, then
the equations we need to solve are
\[
V^k(z^0,\ldots,z^n) =0, \quad k=0, \ldots, n.
\]

\subsection{A reduction of dimension for a balanced ODE}
If we are interested not only in the equilibria but also in
a more detailed study of the transient behaviour of the ODE, then
we can obtain a rough idea (and bounds) by making the assumption
of full symmetry, i.e.,\ we assume that the arrival rates $\alpha^A$
and the initial states $x_0^A$ depend on $A$ only through $|A|$.
Then the trajectory itself $x^A_t$
depends only on the cardinality of $A$ and so we can
reduce the ODE to an $(n+1)$-dimensional one.
Such a symmetrised ODE can yield more detailed information on the time
to reach a small neighbourhood of the equilibrium point.
(Note that our bounds in \S \ref{settle} are very crude.)

\subsection{Non-Poissonian assumptions}
It may be more reasonable in practise to assume that the
time it takes for a chunk to be downloaded or swapped
is a random variable with a heavy-tailed distribution.
This is not captured by our model. Indeed, the interaction times
are not part of our model at all.
A new, more detailed, model should be worked out.

However, a crude capture of this phenomenon is the replacement of
the Poisson processes $\Phi_\zeta$ by more
general point processes, perhaps with heavy-tailed inter-event times.
As long as these processes obey a functional law of large numbers,
we can (by possibly modifying the scaling parameters) rephrase
Theorem \ref{ODEapprox} and repeat the proof in this more general
case.



\appendix

\section{Drift calculation}
\label{driftcalc}

We consider the set of vectors
\begin{multline*}
\NN = \{ -e_F\} \cup\{ e_A:A \subset F\}
\cup\{-e_A+e_{A'}: A \subset A' \subset F\}
\\
\cup\{-e_A-e_B+e_{A'}+e_{B'}: A \subset A' \subset F,
B \subset B' \subset F, A'-A \subset B, B'-B \subset A\}.
\end{multline*}
For each $\zeta \in \NN$ we define a polynomial $Q_\zeta(x)$,
by comparing \eqref{Qdef} and \eqref{rates}:
\begin{align}                        
Q_{e_A}(x) &:= \alpha^A  \nonumber \\
Q_{-e_F}(x) &:= \delta x^F   \nonumber \\
Q_{-e_A+e_{A'}}(x) &:= \lambda_{A,A'}(x)   \nonumber \\
Q_{-e_A-e_B+e_{A'}+e_{B'}}(x) &:= \mu_{A,B}(x) 
\delta_{A,A',B,B'}, 
\label{Qdef2}
\end{align}
where $\lambda_{A,A'}(x)$, $\mu_{A,B}(x)$ are given by \eqref{lll}, 
\eqref{mmm}, respectively, and
\[
\delta_{A,A',B,B'} := 
\1(
A \sqsubset A', 
A'-A \subset B, 
B \sqsubset B', 
B'-B \subset A).
\]
The variable $x$ ranges in $\Z_+^{\PP(F)}$
or in $\R_+^{\PP(F)}$. The algebra is the same in both cases.
Define the drift vector field by $\sum_y (y-x) q(x,y)$. Comparing
\eqref{Qdef2} and \eqref{rates} we have
\[
\sum_y (y-x) q(x,y) = \sum_{\zeta \in \NN} \zeta Q_\zeta(x).
\]
The latter sum appears in \eqref{inteq}, in the course of the proof of Theorem
\ref{ODEapprox}.
We shall verify that
\[
u(x):= \sum_{\zeta \in \NN} \zeta Q_\zeta(x) = v(x),
\]
where $v(x)$ is defined by \eqref{vcomp}.

Consider the terms in the summation $u(x)=\sum_{\zeta \in \NN} \zeta Q_\zeta(x)$ involving $\zeta=-e_A-e_B+e_{A'}+e_{B'}$.
Notice that swapping $A$ with $B$ or $A'$ with $B'$ will not change
the value of $x-e_A,x-e_B+e_{A'}+e_B$, so we need to make sure
to take into account this change only once in the summation.
If we {\em simultaneously} swap $A$ with $B$ {\em and} $A'$ with $B'$
then neither $x-e_A,x-e_B+e_{A'}+e_B$ nor the value
of $Q_{-e_A-e_B+e_{A'}+e_{B'}}(x) = \mu_{A,B}(x) \delta_{A,A',B,B'}$ 
will change because, clearly,
\[
\mu_{A,B}(x) \delta_{A,A',B,B'} 
=
\mu_{B,A}(x) \delta_{B,B',A,A'},
\]
as readily follows from \eqref{lll} and \eqref{mmm}.
We now see that to swap $A$ with $B$ {\em without} swapping $A'$
with $B'$ is impossible (unless $A'=B'$). Indeed, it is
an easy exercise that
\[
\delta_{A,B,A',B'} = \delta_{B,A,A',B'} 
~ \Rightarrow ~ A' = B'.
\]
Taking into account this, we write
\begin{multline}
\label{uuu}
u(x) = \alpha^A e_A -\delta x^F e_F
+ \sum_{A,A'} (-e_A+e_{A'}) \lambda_{A,A'}(x)
\\
+ \frac{1}{2} \sum_{A,B,A',B'} (-e_A-e_B+e_{A'}+e_{B'}) \mu_{A,B}(x)
\delta_{A,B,A',B'},
\end{multline}
where the $1/2$ appears because each term must be counted exactly once.
The variables $A,A',B,B'$ in both summations are  free to move
over $\PP_n(F)$ (but notice that restrictions have effectively
been pushed in the definitions of $\lambda_{A,A'}$, $\mu_{A,B, }$
and $\delta_{A,B,A',B'}$).

Since
\begin{align*}
\sum_{A,B,A',B'} e_A \mu_{A,B}(x) \delta_{A,B,A',B'}
&=
\sum_{A,B,A',B'} e_B \mu_{A,B}(x) \delta_{A,B,A',B'},
\\
\sum_{A,B,A',B'} e_{A'} \mu_{A,B}(x) \delta_{A,B,A',B'}
&=
\sum_{A,B,A',B'} e_{B'} \mu_{A,B}(x) \delta_{A,B,A',B'},
\end{align*}
we have
\begin{multline*}
v(x) = \alpha^A e_A -\delta x^F e_F -\sum_{A,A'} e_A \lambda_{A,A'}(x)
+ \sum_{A,A'} e_{A'} \lambda_{A,A'}(x)
\\
- \sum_{A,B,A',B'} e_A \mu_{A,B}(x) \delta_{A,B,A',B'}
+ \sum_{A,B,A',B'} e_{A'} \mu_{A,B}(x) \delta_{A,B,A',B'}.
\end{multline*}
Call the four sums appearing in this display as
$u\I(x) , u\II(x) , u\III(x) , u\IV(x)$, in this order.
We use the definitions \eqref{lll}, \eqref{mmm} of $\lambda_{A,A'}$, 
$\mu_{A,B}$
and find the components of the vectors $u\I, \ldots, u\IV$
by hitting each one with a unit vector $e_G$,
i.e.\ by taking the inner products
$v^G\I = \langle e_G, u\I\rangle, \ldots, v^G\IV = \langle e_G, u\IV\rangle$.
We have:
\begin{align}
u\I^G(x) 
= 
-\sum_{A'} \lambda_{G,A'}(x)
&=
-\sum_{A'} \beta x^G
\sum_{B : B \supset A'}
\frac{x^B}{|B-G|} \1(G \sqsubset A')
                                                \nonumber
\\
&=
-\beta x^G
\sum_B \frac{x^B}{|B-G|} \sum_{A'} \1(G \sqsubset A' \subset B)
=
-\beta x^G
\sum_{B \supset G} x^B ,
                                        \label{I}
\end{align}
where, in deriving the last equality we just observed that the number of
sets $A'$ that contain one more element than $G$ and are contained
in $B$ is equal to $|B-G|$, as long as $G \subset B$:
\[
\sum_{A'} \1(G \sqsubset A' \subset B) = |B-G| ~\1(G \subset B).
\]
Next,
\begin{align*}
u\II^G(x) 
= \sum_A \lambda_{A,G}(x)
&=
\sum_A
\beta x^A
\sum_{B \supset G}
\frac{x^B}{|B-A|}
~
\1(A \sqsubset G)
\\
&= \beta \sum_{B \supset G} x^B \sum_A \frac{x^A}{|B-A|}
~\1(A \sqsubset G)
\end{align*}
Notice that, in the last summation, $G$ contains exactly one more
element than $A$ and is strictly contained in $B$, so $|B-A|=|B-G|+1$.
Hence
\begin{equation}
u\II^G(x)
=
\beta \sum_{B : B \supset G} \frac{x^B}{|B-G|+1} 
\sum_A x^A \1(A \sqsubset G)
= \beta \sum_{B : B \supset G} \frac{x^B}{|B-G|+1}
\sum_{g \in G} x^{G-g}.
                                        \label{II}
\end{equation}
For $u\III(x)$, we have:
\begin{multline}
u\III^G(x) 
=
- \sum_{B,A',B'} \mu_{G,B}(x) \delta_{G,A',B,B'}
\\
= 
-\gamma 
\sum_B \frac{x^G x^B}{|G \setminus B| |B \setminus G|}
\cdot \sum_{A'} \1(G \sqsubset A', A'-G \subset B)
\cdot \sum_{B'} \1(B \sqsubset B', B'-B \subset G)
\\
=
-\gamma
\sum_B \frac{x^G x^B}{|G \setminus B| |B \setminus G|}
\cdot |B \setminus G|~ \1(B \setminus G \not = \varnothing)
\cdot |G \setminus B|~ \1(G \setminus B \not = \varnothing)
\\
= -\gamma x^G \sum_{B \not \sim G} x^B.
                        \label{III}
\end{multline}
As for the last term, we have:
\begin{align}
u\IV^G(x) 
&= 
\sum_{A,B,B'} \mu_{A,B}(x) \delta_{A,G,B,B'}
\nonumber
\\
&=
\gamma \sum_B \sum_A
\frac{x^A x^B}{|A \setminus B| |B \setminus A|}
\1(A \sqsubset G, G-A \subset B)
\sum_{B'} \1(B \sqsubset B', B'-B \subset A)
\nonumber
\\
&=
\gamma  \sum_B \sum_A
\frac{x^A x^B}{|B \setminus A|}
\1(A \sqsubset G, G-A \subset B)
~
\1(A \setminus B \not = \varnothing)
\nonumber
\\
&=
\gamma \sum_B \frac{x^B}{|B \setminus G|+1}
\sum_A x^A
~\1(A \sqsubset G, G-A \subset B)
~\1(G \not \subset B)
\nonumber
\\
&=
\gamma \sum_B \frac{x^B}{|B \setminus G|+1}
\sum_{g \in G \cap B} x^{G-g}
~\1(G \not \subset B)
                                        \label{IV}
\end{align}
Adding \eqref{I} and \eqref{III} we obtain the first part of
\eqref{vcomp}, while \eqref{III} and \eqref{IV} give the
second part.

\end{document}